\documentclass[11pt]{article}
\usepackage{amsmath,amsthm,amsfonts,verbatim,amssymb}

\theoremstyle{plain}
 \newtheorem{lemma}{Lemma}[section]
 \newtheorem{theorem}[lemma]{Theorem}
 \newtheorem{corollary}[lemma]{Corollary}

 \theoremstyle{definition}

\begin{document}

\title{Embedding the diamond graph in $L_p$ and \\ dimension reduction in $L_1$}

\author{
    James R. Lee\thanks{Work partially supported by NSF grant CCR-0121555
    and an NSF Graduate Research Fellowship.} \\
    U.C. Berkeley
    \and
    Assaf Naor \\
    Microsoft Research
     }

\date{\empty}

\maketitle

\begin{abstract}
We show that any embedding of the level $k$ diamond graph of
Newman and Rabinovich~\cite{nr} into $L_p$, $1 < p \le 2$,
requires distortion at least $\sqrt{k(p-1) + 1}$.  An immediate
corollary is that there exist arbitrarily large $n$-point sets $X
\subseteq L_1$ such that any $D$-embedding of $X$ into $\ell_1^d$
requires $d \geq n^{\Omega(1/D^2)}$. This gives a simple proof of
a recent result of Brinkman and Charikar~\cite{bc} which settles
the long standing question of whether there is an $L_1$ analogue
of the Johnson-Lindenstrauss dimension reduction lemma~\cite{jl}.

\end{abstract}

\section{The diamond graphs, distortion, and dimension}

We recall the definition of the diamond graphs
$\{G_k\}_{k=0}^\infty$ whose shortest path metrics are known to be
uniformly bi-lipschitz equivalent to a subset of $L_1$
(see~\cite{sinclair} for the $L_1$ embeddability of general series-parallel graphs). The diamond graphs were used in~\cite{nr} to
obtain lower bounds for the Euclidean distortion of planar graphs
 and similar arguments were previously used in a different context by Laakso~\cite{laakso}.

$G_0$ consists of a single edge of length 1.  $G_{i}$ is obtained
from $G_{i-1}$ as follows. Given an edge $(u,v) \in E(G_{i-1})$,
it is replaced by a quadrilateral $u,a,v,b$ with edge lengths
$2^{-i}$. In what follows, $(u,v)$ is called an edge of level
$i-1$, and $(a,b)$ is called the level $i$ anti-edge corresponding
to $(u,v)$.  Our main result is a lower bound on the distortion
necessary to embed $G_k$ into $L_p$, for $1 < p \leq 2$.

\begin{theorem}\label{thm:ONLY} For every $1<p\le 2$,
any embedding of $G_k$ into $L_p$ incurs distortion at least $
\sqrt{1+(p-1)k}$.
\end{theorem}

The following corollary shows that the diamond graphs cannot be
well-embedded into low-dimensional $\ell_1$ spaces.
In particular,
an $L_1$ analogue of the Johnson-Lindenstrauss dimension reduction
lemma does not exist.
The same graphs were used in~\cite{bc} as an example
which shows the impossibility of dimension reduction in $L_1$. Our
proof is different and, unlike the linear programming
based argument appearing there, relies on geometric intuition.
We proceed by observing that a lower bound
on the rate of decay of the distortion
as $p \to 1$ yields a lower bound on the required dimension in $\ell_1$.

\begin{corollary}
For every $n \in \mathbb N$, there exists an $n$-point subset $X
\subseteq L_1$ such that for every $D > 1$, if $X$ $D$-embeds into
$\ell_1^d$, then $d \geq n^{\Omega(1/D^2)}$.
\end{corollary}

\begin{proof}
Since $\ell_1^d$ is $O(1)$-isomorphic to $\ell_{p}^d$ when $p = 1
+ \frac{1}{\log d}$ and $G_k$ is $O(1)$-equivalent to a subset $X
\subseteq L_1$, it follows that $\sqrt{1+ \frac{k}{\log d}}
= O(D)$.  Noting that $k = \Omega(\log n)$ completes the proof.
\end{proof}

\section{Proof}

The proof is based on the following inequality.
The case $p=2$ is the well known ``short diagonals lemma'' which was
central to the argument in~\cite{laakso,nr}.

\begin{lemma}\label{lemma:us}
 Fix $1<p\le 2$ and $x,y,z,w\in L_p$. Then,
$$
\|y-z\|_p^2+(p-1)\|x-w\|_p^2\le
\|x-y\|_p^2+\|y-w\|_p^2+\|w-z\|_p^2+\|z-x\|_p^2.
$$
\end{lemma}

\begin{proof} For every $a,b\in
L_p$, $||a+b||_p^2 + (p-1) ||a-b||_p^2 \leq 2 (||a||_p^2 +
||b||^2_p)$. A simple proof of this classical fact can be found,
for example, in~\cite{bcl}. Now,
$$||y-z||_p^2 + (p-1) ||y - 2x +
z||_p^2 \leq 2 ||y-x||_p^2 + 2 ||x-z||_p^2$$ and $$||y-z||_p^2 +
(p-1) ||y - 2w + z||_p^2 \leq 2 ||y-w||_p^2 + 2||w-z||_p^2.$$
Averaging these two inequalities yields
$$||y-z||_p^2 + (p-1)\frac{||y -2x + z||_p^2 + ||y -2w+
z||_p^2}{2} \leq
\|x-y\|_p^2+\|y-w\|_p^2+\|w-z\|_p^2+\|z-x\|_p^2.$$ The required
inequality follows by convexity. 
\end{proof}

\begin{lemma}\label{lem:sum}
Let $A_i$ denote the set of anti-edges at level $i$ and let
$\{s,t\} = V(G_0)$, then for $1<p \leq 2$ and any $f: G_k \to
L_p$,
$$
||f(s) - f(t)||_p^2 + (p-1)\sum_{i=1}^k \sum_{(x,y) \in A_i}
||f(x) - f(y)||_p^2 \leq \sum_{(x,y) \in E(G_k)} ||f(x) -
f(y)||_p^2.$$
\end{lemma}

\begin{proof}
Let $(a,b)$ be an edge of level $i$ and $(c,d)$ its corresponding
anti-edge.  By Lemma \ref{lemma:us}, $||f(a)-f(b)||_p^2 +
(p-1)||f(c)-f(d)||_p^2 \leq ||f(a)-f(c)||_p^2 + ||f(b)-f(c)||_p^2 + ||f(d)-f(a)||_p^2 +
||f(d)-f(b)||_p^2$.  Summing over all such edges and all $i = 0, \ldots,
k-1$ yields the desired result by noting that the terms $||f(x) - f(y)||_p^2$
corresponding to $(x,y) \in E(G_i)$ cancel for $i = 1, \ldots, k-1$.
\end{proof}

\noindent
The main theorem now follows easily.

\begin{proof}[Proof of Theorem~\ref{thm:ONLY}]
Let $f:G_k \to L_p$ be a non-expansive $D$-embedding.  Since
$|A_i| = 4^{i-1}$ and the length of a level $i$ anti-edge is
$2^{1-i}$, applying Lemma \ref{lem:sum} yields
$\frac{1+(p-1)k}{D^2} \leq 1$.
\end{proof}

\end{document}